\def\nat{{\tt I\kern-.2em{N}}}
\def\hyper#1{\ ^*\kern-.2em{#1}}

\def\lint{\lceil}
\def\rint{\rceil}
\def\qed{{\vrule height6pt width3pt depth2pt}\par\medskip}
\def\sig{{^\sigma}}
\def\parm{\par\medskip}
\def\pars{\par\smallskip}
\def\hypernat{{^*{\nat }}}
\def\m@th{\mathsurround=0pt}

\magnification=\magstep1
\tolerance 10000
\font\eightrm=cmr9
\baselineskip  12pt
\centerline{\bf Modeling the Dialectic}\par\bigskip 
\centerline{Robert A. Herrmann}\par\medskip
\centerline{26 OCT 2008}\par\bigskip
{\leftskip=0.5in \rightskip=0.5in \noindent {\eightrm {\it Abstract:} Models are constructed that satisfy each   
dialectical scheme $TAS_1,\ TAS_2,\ TAS_3.$ Significantly different finite models that satisfy $TAS_1$ and $TAS_2$, a denumerable model that satisfies schemes $TAS_i,\ i=1,2,3$ and an infinite hyperfinite model, with a single antithesis, that satisfies $TAS_1$ and $TAS_2$ are defined. It is shown that no finite model satisfies $TAS_3$.  \par}}\pars 
Mathematics Subject Classifications 03B22, 03B65.\parm

\noindent {\bf 1. Introduction.} \parm 

Three dialectical schemes can be formally expressed in a 
first-order language with equality (Gagnon, 1980). Consider
a set of predicates $T(-),\ N(-),$ $A(-,-),\ D(-,-),$ $P(-
,-),\ S(-,-,-).$  The three formal schemes are:\parm
\centerline{$TAS_1$}
\line{\indent\indent E1 $\exists x[T(x)].$\hfil}
\line{\indent\indent E2 $\forall x[T(x) \to \exists !\, y[A(y,x)]].$\hfil}
\line{\indent\indent E3 $\forall x\forall y[A(y,x) \to \exists !\, z[S(z,x,y)]].$\hfil}
\line{\indent\indent R1 $\forall x\forall y[A(y,x) \to [T(x) \land \neg 
A(x,y)]]$\hfil}
\line{\indent\indent R2 $\forall x\forall y\forall z[S(z,x,y) \to [T(z) \land 
[A(x,y) \lor A(y,x)]\land \neg  [S(x,z,y) \lor S(y,x,z)]]].$\hfil} \pars

\centerline{$TAS_2$}
\line{\indent\indent $N(z){\buildrel \rm def \over =} \exists x\exists 
y[S(z,x,y) \land D(z,x) \land D(z,y)].$ \hfil}
\line{\indent\indent E1 $\exists x[T(x)].$\hfil}
\line{\indent\indent E2 $\forall x[T(x) \to \exists !\, y[A(y,x)]].$\hfil}
\line{\indent\indent E3 $\forall x\forall y[A(y,x) \to \exists !\, z[S(z,x,y)]].$\hfil}
\line{\indent\indent E4 $\exists x[N(x)].$\hfil}
\line{\indent\indent E5 $\forall x[N(x) \to \exists y[N(y)\land y \not= x]].$\hfil}
\line{\indent\indent R1 $\forall x\forall y[A(y,x) \to [T(x) \land \neg 
A(x,y)]]$\hfil}
\line{\indent\indent R2 $\forall x\forall y\forall z[S(z,x,y) \to [T(z) \land 
[A(x,y) \lor A(y,x)]\land \neg  [S(x,z,y) \lor S(y,x,z)]]].$\hfil}\pars

\centerline{$TAS_3$}
\line{\indent\indent $N(z){\buildrel \rm def \over =} \exists x\exists 
y[S(z,x,y) \land D(z,x) \land D(z,y)].$ \hfil}
\line{\indent\indent E1 $\exists x[T(x)].$\hfil}
\line{\indent\indent E2 $\forall x[T(x) \to \exists !\, y[A(y,x)]].$\hfil}
\line{\indent\indent E3 $\forall x\forall y[A(y,x) \to \exists !\, z[S(z,x,y)]].$\hfil}
\line{\indent\indent E4 $\exists x[N(x)].$\hfil}
\line{\indent\indent E5.1 $\forall x[N(x) \to \exists y[N(y)\land P(x,y)]].$\hfil}
\line{\indent\indent E6.1 $\forall x\exists y [P(x,y)].$\hfil}
\line{\indent\indent R1.1 $\forall x\forall y[A(y,x) \to [T(x) \land P(x,y)]].$\hfil}
\line{\indent\indent R2.1 $\forall x\forall y\forall z[S(z,x,y) \to [T(z) \land 
[A(x,y) \lor A(y,x)] \land S(z,y,x) \land P(x,z) \land P(y,z)]].$\hfil}
\line{\indent\indent R3.1 $\forall x\forall y[P(x,y) \to \neg P(y,x)].$\hfil}
\line{\indent\indent R4.1 $\forall x\forall y\forall z[[P(x,y)\land P(y,z) ]\to 
P(x,z)]].$\hfil} \pars
\parm\vfil\eject
 
\noindent {\bf 2. Standard Models.} \parm
For any model, the axioms require all relations to be nonempty. Let nonempty $T$ be the set of theses and nonempty ${\cal A}$ be the set of antitheses. In this section, each of the designated models has domain $T \cup {\cal A}.$ For the dialectic, $T(a)$ is interpreted ($\lint T(a)\rint$): $a$ is a theses. Further, $\lint A(b,a)\rint$: $b$ is the antitheses of $a$, $\lint S(a,b,c)\rint$: $a$ is the synthesis of $b$ and $c$, $\lint D(a,b)\rint$: $a$ is qualitatively different than $b$ and $\lint P(a,b) \rint$: represents an order for $a$ and $b.$ This order is often related to ``time.'' The defined predicate $N$ restricted to various theses yields a thesis termed a ``nodel point.'' The axioms imply that if $TAS_1$ and $TAS_2$ have finite models, then R2 requires that each domain contain three or more elements. In what follows, the constant predicate symbol is used for the corresponding set theoretic object. \parm

\centerline{\bf Model A}\parm

\noindent{\bf Definition 2.1.} The numbers 1, 2, 3 are considered as but distinct symbols. The ``='' means identical as symbols.\pars
(a) Let $T^A =\{1,2,3\}= {\cal A}^A,$\pars
(b) $A^A = \{(1,2),(2,3),(3,1)\},\ S^A=\{(1,3,2),(2,1,3),(3,2,1)\},\ D^A = \hfil\break
 \{(1,3),(1,2), (2,1), (2,3),(3,2),(3,1)\}.$\pars
(c) $N^A = \{1,2,3\}$ (or $\{1,2\}).$\parm

\noindent {\bf Theorem 2.2.} {\it The structure ${\cal D}^A = \langle T^A, S^A, A^A, D^A,N^A \rangle$ is a model for $TAS_1$ and $TAS_2.$ Hence, $TAS_1$ and $TAS_2$ are, relative to models, consistent schemes.} \pars
Proof. For E1, $1 \in T^A.$ For E2, let $1 \in T^A.$ Then $(x,1) \in A^A$ if and only if $x =3.$ 
Let $2 \in T^A.$ Then $(x,2) \in A^A$ if and only if $x = 1.$ Let $3 \in T^A$. Then $(x,3) \in A^A$ if and only if $x = 2.$ For E3, let $(1,2) \in A^A.$ Then $(x,2,1)\in S^A$ if and only if $x = 3.$ Let $(2,3) \in A^A.$ Then $(x,3,2) \in S^A$ if and only if $x = 1.$ Let $(3,1) \in A^A.$ Then
$(x, 1,3) \in S^A$ if and only if $x = 2.$ R1 is obvious. For R2, let $(3,2,1) \in S^A.$ Then $3 \in T^A,\ (1,2) \in A^A, (2,3,1)\notin S^A,\ (1,2,3) \notin S^A.$ Let $(1,3,2) \in S^A.$ Then $1 \in T^A,\ (2,3) \in A^A,\ (3,1,2)\notin S^A,\ (2,3,1) \notin S^A.$ Let $(2,1,3) \in S^A.$ Then $2 \in T^A,\ (3,1) \in A^A,\ (1,2,3) \notin S^A,\ (3,1,2) \notin S.$ Hence (a) (b) model $TAS_1.$ \pars

For E4, $1 \in N^A.$ E5 is obvious for $N^A$ (or $\{1,2\}$), since $1 \not= 2 \in N^A.$ Hence, (a) (b) (c) model $TAS_2$. \qed \parm

\centerline{\bf Model B.}\parm

\noindent{\bf Definition 2.3.} The numbers 1, 2, 3, 4 are considered as but distinct symbols. The ``='' means identical as symbols. \pars
(a) Let $T^B \{1,2\},\ {\cal A}^B = \{3,4\},$\pars
(b) $A^B = \{(3,1),(4,2)\},\ S^B =\{(2,1,3),(1,2,4)\},\ D^B= \{(2,1),(2,3), (1,2), (1,4)\},$\pars
(c) $N^B = \{1,2\}$\parm

\noindent {\bf Theorem 2.4.} {\it The structure ${\cal D}^B = \langle T^B \cup {\cal A}^B, S^B, A^B, D^B,N^B \rangle$ is a model for $TAS_1$ and $TAS_2.$} \pars
Proof. For E1, $1 \in T^B.$ For E2, let $1 \in T^B.$ Then $(x,1) \in A^B$ if and only if $x =3.$ 
Let $2 \in T^B.$ Then $(x,2) \in A^B$ if and only if $x = 4.$ For E3, let $(3,1) \in A^B.$ Then $(x,1,3)\in S^B$ if and only if $x = 2.$ Let $(4,2) \in A^B.$ Then $(x,2,4) \in S^B$ if and only if $x = 1.$ 
R1 is obvious. For R2, let $(2,1,3) \in S^B.$ Then $2 \in T^B,\ (3,1) \in A^B, (1,2,3)\notin S^B,\ (3,1,2) \notin S^B.$ Let $(1,2,4) \in S^B.$ Then $1 \in T^B,\ (4,2) \in A^B,\ (2,1,4)\notin S^B,\ (4,2,1) \notin S^B.$ Hence (a) (b) model $TAS_1.$ \pars

For E4, $1 \in N^B.$ For E5, $1 \not= 2 \in N^B.$ Hence, (a) (b) (c) model $TAS_2$. \qed \parm

\centerline{\bf Model C}\parm

\noindent {\bf Definition 2.5.} For the natural numbers $\nat$, let $a_0= 3,\ b_0 = 4,\ c_0 = a_0 + b_0 = 7.$ By induction over $\nat$, define\pars
 
(a) $\forall i \geq 0,\ a_{i+1} = c_i,$\pars
(b) $\forall i \geq 0,\ b_{i+1} = c_i + 1 (= a_{i+1} + 1).$\pars
(c) $\forall i \geq 0,\ c_{i+1} = 2c_i + 1.$\parm

From (c) and (b), $\forall i \geq 0,\ c_{i + 1} = c_i + c_i + 1 = a_{i+1} + b_{i+1}.$ From the initial values, this yields that (1) $\forall i \geq 0, c_i = a_i + b_i.$ Further, $\forall i\geq 0,\ b_i = a_i +1 > a_i.$ Hence, (2) $\forall i \geq 0,\ a_{i+1} = c_i = a_i + b_i = a_i + a_i + 1 = 2a_i + 1.$ Then (3) $\forall i \geq 0,\ b_{i+1} = a_{i+1} + 1 = 2a_i + 1 + 1 = 2a_i +2,$ and $c_i > a_i,\ c_i > b_i.$ \parm

\noindent{\bf Definition 2.6.} For the theses, let 
$T^C=\{a_i\mid i \in \nat\}$; for the 
antitheses, let ${\cal A}^C=\{b_i\mid i \in \nat\}.$ (Then $T^C \cap {\cal A}^C = \emptyset.$) Each of the remaining undefined predicates is modeled by an appropriate relation.\parm
(1) Let $S^C = \{(c_i,a_i,b_i)\mid i \in \nat\} \cup \{(c_i,b_i,a_i)\mid i \in 
\nat\},$\pars
(2) $A^C = \{(b_i,a_i)\mid i \in \nat\},$\pars
(3) $D^C = \{(c_i,b_i)\mid i \in \nat\} \cup \{(c_i,a_i)\mid i \in \nat\}.$ Then $N^C = \{c_i\mid i \in \nat\} = \{a_{i+1}\mid i \in \nat\},$\pars
(4) $P^C = \{(x,y)\mid [x, y \in  T \cup {\cal A}] \land [x<y] \}.$\parm

\noindent{\bf Theorem 2.7.} {\it The structure ${\cal D}^C=\langle T^C \cup {\cal A}^C,T^C,S^C,A^C,D^C,N^C,P^C\rangle$ 
is a model for $TAS_i,\ i = 1,2,3.$}\pars
Proof.   E1, E2, E3, E4, and R1 are obvious. Note that $\forall i \geq 0, c_i > a_i,\ c_i > b_i.$ Hence, $\neg [(x,z,y)\in S^C \lor (y,x,z)\in S^C]$ holds for $T^C \cup {\cal A}^C.$ From this, R2 holds. For each $\forall i \geq 0$ and  $c_i,$ there exists a $c_{i + 1} > c_i$ and $c_{i+1} \not= c_i.$ Hence, E5 and E5.1 hold. Since $\forall i \geq 0, \ a_i < b_i,$ R1.1 holds. Since $\forall i \geq 0,\ a_i < c_i,\ b_i < c_i$, then this and definition 2.6 part (1) imply that R2.1 holds and E6.1 holds from the definition of $P^C$. From the properties of the order $<$, R3.1 and R4.1 hold. \qed 

Theorem 2.7 implies that $TAS_3$ is consistent relative to the theory of natural numbers. Gagnon (1980) uses the theory of $\nat$ and states that $TAS_i,\ i = 1,2,3$ satisfy a different denumerable model. Moreover, he states that other dialectical theories can be generated from $TAS_3$ by adding axioms. For a set X, the term ``finite'' means that either $X = \emptyset$ or for some $n \in \nat, n \geq 1,$ there exists a bijection $f\colon [1,n] \to X.$ Gagnon does not show that $TAS_1,\ TAS_2$ have finite models. Gagnon does not mention that scheme $TAS_3$ has no finite model. For a nonempty set $X$ and $1 \leq n \in \nat$, the notation $\vert X \vert = n$ signifies that there exists a bijection $f\colon [1,n] \to X.$ \parm 

\noindent{\bf Theorem 2.8} {\it There does not exist a finite model that satisfies scheme $TAS_3$.}\pars

Proof. It should be well known that if a set $N$ satisfies E4, E5.1 and binary relation $P$ satisfies R3.1 and R4.1, then $N$ is (ordinary) infinite.  Indeed, $N$ is nonempty by E4. Let $p \in N$ and $(1,p) \in f_3.$ Then there exists $q \in N$ such that $(p,q) \in P.$ Then from R3.1, $p \not= q.$ Let $(2,q) \in f_3.$ Then there exists $r\in N$ such that $(q,r) \in P.$ and $q \not= r.$ If $r= p,$ then R4.1 implies that $(p,p)\in P;$ a contradiction. Let $(3,r)\in f_3$. (1) For natural number $n = 3,$ there exists an injection $f_3\colon [1,3] \to N$ such that if $i,j \in \nat$, where $1\leq i,j \leq 3,\ i \not= j,$ then $f_3(i) \not= f_3(j)$. If $1\leq i < j\leq 3,$ then $(f_3(i),f_3(j))\in P.$ Hence, $\vert N \vert \geq 3$ or $N$ is infinite.\pars

(2) Assume that for $n \geq 3,$ there exists an injection $f_n \colon [1,n] \to N$ such that if $i,j \in \nat$, where $1\leq i,j \leq n,\ i \not= j,$ then $f_n(i) \not= f_n(j)$. If $1\leq i < j\leq n,$ then $(f_n(i),f_n(j))\in P.$  (3) For $n\geq 3,$ consider $n +1.$ Then there exists an injection $f_n \colon [1,n] \to N$ such that if $i,j \in \nat$, where $1\leq i,j \leq n,\ i \not= j,$ then $f_n(i) \not= f_n(j)$. If $1\leq i < j\leq n,$ then $(f_n(i),f_n(j))\in P.$ Since $f_n(n) \in N$, then there exists $s \in N$ such that $(f_n(n), s)\in P.$ If $s =f_n(n),$ R3.1 is contradicted. Hence, $f_n(n) \not= s.$ If $s = f_n(j), 1\leq j < n,$ then $(s, f_n(n))\in P.$ Hence, $(s,s) \in P;$ a contradiction. Let $f_{n+1} = f_n \cup \{(n+1,s)\}.$ Then $f_{n+1}$ is an injection on $[1,n+1]$ into $N$ and it follows from (2) that if $i,j \in \nat$, where $1\leq i,j \leq n+1,\ i \not= j,$ then $f_{n+1}(i)\not= f_{n+1}(j)$and if $1\leq i < j\leq n+1,$ then $(f_{n+1}(i),f_{n+1}(j))\in P.$ Consequently, by induction, for each $n \geq 3,\ \vert N \vert \geq n$ or $N$ is infinite. If $1 \leq \vert N \vert = m \in \nat$, then $\vert N \vert \geq m+2;$ a contradiction since finite cardinalities satisfy natural number order properties. Hence, $N$ is infinite. \qed

\centerline{\bf Model D}\parm

In what follows, infinitely many finite models for $TAS_1$ and $TAS_2$ are defined and each has the special property that there is but one antithesis. Simple properties of the natural numbers are used. \parm

\noindent{\bf Definition 2.9} Let $k \in \nat,\ k \geq 2.$\pars 
(a) Let $T_k^D = \{i\mid (i \in \nat)\land(0\leq i \leq k)\}$, ${\cal A}_k^D = \{k+1\},$\pars
(b) $A_k^D = \{(k+1,i)\mid (i\in \nat)\land(0\leq i \leq k)\},\ S_k^D = S'_k \cup S''_k,$ where $S'_k = \{(i+1, i, k+1)\mid (i \in \nat)\land(0 \leq i \leq k-1)\}$ and $S''_k = \{(0,k,k+1)\},$ $D_k^D = \{(1,0),(1,k+1),(2,1),(2,k+1)\},$\pars
(c) $ N^D = \{1,2\}.$\parm
\noindent{\bf Theorem 2.10} {\it For each $k \in \nat,\ k \geq 2,$ the structure ${\cal D}_k^D=\langle T_k^D\cup {\cal A}_k^D, T_k^D,{\cal A}_k^D, S_k^D, A_k^D,\break N^D,D_k^D \rangle$ is a model for $TAS_1$ and $TAS_2.$}\pars
Proof. (Although the defined ${\cal A}_k^D$ need not be included in ${\cal D}_k^D$, it is useful to list it.) For E1, $0 \in T_k^D$. For E2, let $i\in T_k^D,\ 0\leq i\leq k.$ Then $(x,i) \in A_k^D$ if and only if $x = k+1.$ For E3, let $i = k,\ (k+1,k) \in A_k^D.$ Then $(x,k,k+1)\in S_k^D$ if and only if $x = 0.$ Let $0\leq i \leq k-1,\ (k+1,i) \in A_k^D.$  Then $(x, i, k+1) \in S_k^D$ if and only if $x = i+1.$ R1 holds since for $ 0\leq i \leq k,\ (i,k+1)\notin A_k^D.$ For R2, let $(0,k,k+1) \in S_k^D.$ Then $0 \in T_k^D, (k+1,k)\in A_k^D$ and, since $k >1,$ $(k,0,k+1) \notin S_k^D,\ (k+1,k,0)\notin S_k^D.$ Let $(i+1,i, k+1) \in S_k^D,\ 0\leq i\leq k-1.$ Then $i+1 \in T_k^D,\  (k+1,i) \in A_k^D$ and $(i,i+1,k+1) \notin S_k^D,\ (k+1,i,i+1) \notin S_k^D.$ Hence, (a) and (b) model $TAS_1.$\pars
For E4, $N^D \subset T_k^D$ and $\vert N^D\vert = 2.$ Moreover, the 1 and 2 are not just distinct as natural numbers but have additional qualities. The 1 is an odd number and a multiplicative identity. The 2 is the successor, even and not an identity. \qed

The fact that $\cal A$ can contain but one member does not hold for a model for $TAS_3.$ Consider a model for $TAS_3.$ Let $E$ be the domain. Then $\emptyset \not= T \subset E.$ Let ${\cal A} = \{y\mid (\exists x \in T)\land((y,x) \in A)\}.$ Then ${\cal A} \subset E.$ Suppose that $B \subset {\cal A}$ and $\vert {B} \vert = 1.$ Let $b \in {B}.$ Then there exists an $x \in T$ such that $(b,x) \in A.$ By R1.1, $(x,b) \in P$ and by R3.1 $x \not= b.$ By E3, there exists a $z \in T$ such that $(z,x,b) \in S.$ Since R2.1 holds for $S$, then $(x,z), \ (b,z) \in P.$ Hence, $z\not= b,\ z\not= x$ by R3.1. There exists a $c \in {\cal A}$ such that $(c,z) \in A, \ c\not= z.$ If $c \in B$, then $c = b$ and $(z,b) \in P$ by R1.1. By R4.1 $(b,b)\in P.$ This contradicts R3.1. Hence, $\vert {\cal A}\vert \geq 2$ or is infinite.\pars

For any model for $TAS_3,$ the domain for the binary relation $P^C$ needs to be specified in order to apply E5.1. From E5.1, $P^C \cap N^C \times N^C \not= \emptyset.$ From E2 and R1.1, $P^C \cap  T^C \times {\cal A}^C\not= \emptyset.$ Define ${\cal A}(x) = \exists y[y \in T\land A(x,y)].$ Obviously, if the axiom $\forall x[{\cal A}(x) \to \exists y [{\cal A}(y) \land P(x,y)]]$ is added to $TAS_3,$ then $P^C \cap {\cal A}^C \times {\cal A}^C \not= \emptyset$ and ${\cal A}^C$ is infinite. 
  \parm

\noindent{\bf 3. Hyperfinite Models.} \parm

Let $\hyper {\cal M} = \langle \kern-.3em\hyper {X},\in,= \rangle$ be a Robinson-styled nonstandard model for all bounded set theoretic first-order statements that hold in $\langle {X}, \in,=\rangle$, where $X$ is a superstructure with atoms $A$ and $\nat \subset A$
(Herrmann, 1991, 1993). The following result describes a special collection of ultradialectics. \parm
\noindent{\bf Theorem 3.1} {\it Each infinite natural number $\lambda$ generates an infinite hyperfinite structure ${\cal D}_\lambda^H =\langle T_\lambda^H \cup {\cal A}_\lambda^H, T_\lambda^H,{\cal A}_\lambda^H = \{\lambda + 1\},\kern-.7em\hyper N^H,S_\lambda^H, A_\lambda^H ,D_\lambda^H \rangle$. Each $Y \in \{T_\lambda^H \cup {\cal A}_\lambda^H,T_\lambda^H,{\cal A}_\lambda^H, \break \kern-.4em\hyper N^H, S_\lambda^H, A_\lambda^H,D_\lambda^H\}$ is an internal member of \kern-.2em$\hyper {X}$, and ${\cal D}_\lambda^H$ models $TAS_1$ and $TAS_2$.}\pars
Proof. Simply let $\lambda \in \hypernat - \sig\nat$, or $\hypernat - \nat$ when each $\hyper p \in \sig\nat$ is identified with $p \in \nat.$ In definitions 
(a) (b) in 2.9, replace $k$ with $\lambda$. Since the natural number properties used in the proof that establishes Theorem 2.10 hold for members of $\hypernat$ (Herrmann 1991), then all the axioms for $TAS_1$ and $TAS_2$ hold for ${\cal D}_\lambda^H$. Each $Y$ is equal to a set expressed using the internal definition principle (Herrmann, 1991, p. 29). For example, $S'_\lambda = \{(z,x,y)\mid (z \in \hypernat)\land(x \in \hypernat) \land (y \in \hypernat)\land (0 \leq x \leq \lambda-1)\land(z = x +1)\land (y = \lambda +1)\}.$ \qed 

\noindent{\bf 4. General Logic-Systems.}\parm

In any model for these three schemes, the $S$ relation is the basic one that models the results of the actual ``logical argument.'' Consider the inverse relation $S^{-1} = \{(y,x,z)\mid (z,x,y) \in S\}.$ The set $\{S^{-1}\}= RI(T\cup {\cal A})$ represents a set of ``finitary rules of inference'' for a finite logic-system defined on the language $T \cup {\cal A}$ (Herrmann, 2001, 2006a). For finite logic-systems, the informal algorithm used to ``deduce'' from the set of premises $X \subset T \cup {\cal A}$ is modeled after the formalizable predict-logic processes used throughout scientific discourse. These processes generate an equivalent ``finite consequence operator.'' Although mostly not displaying the usual forms for predicate deduction, the mental processes employed to produce a formal dialectical argument using any of these three schemes and their extensions are equivalent to those mental processes used for formal predicate deduction (Herrmann, 2001, 2006b). The $S^{-1}$ yields a theory and behavioral-signature (Herrmann, 2006b). Such predicate deduction is said to be ``axiomless.'' Significantly, physical first-order theories are equivalent to theories produced by axiomless first-order predicate deduction (Herrmann,2006b).\par\vfil\eject

\noindent{\bf 5. Applications.}\parm

For an actual denumerable human language $L,$ Gagnon assumes that these schemes apply to events that are both mental and physical. For some events and scheme $TAS_3$, the set $M$ of meaningful stings of symbols is a finite set. Due to Theorem 2.8, a domain $E$ for scheme $TAS_3$ requires denumerably many members of $L.$ In this case, the actual presented dialectical argument would not be applied for any members of $S$ containing any coordinates from $E - M$. However, considering $TAS_3$ applied to mental activity such as mathematics, the set $M$ can be considered as denumerable and a denumerable $E \subset M.$  \pars
The models used for $TAS_1,\ TAS_2$ can be applied to a meaningful language $M$ where each member of $M$ is considered as having distinctly different meanings than any other member. If any of the schemes $TAS_i,\ i=1,2,3,$ is meaningful for a denumerable $E$, then, using the methods in Herrmann (1993), there exists a corresponding distinct ultradialectic.  \pars

Let $f^X\colon (T^X \cup {\cal A}^X) \to M$ be any injection. Then the structure $f^X({\cal D}^X) = \langle f^X[T^X] \cup f^X[{\cal A}^X],f^X[T^X],f^X[S^X],\ldots \rangle,$ where $X = A,B,C,D,$ is (model-theoretic) isomorphic to ${\cal D}^X.$  Hence, each $f^X({\cal D}^X),\ X = A,B,C,D,$ models the corresponding first-order scheme $TAS_i, \ i = 1,2,3.$ [Recall that notation such as $f^X[S^X] = \{(f^X(z),f^X(x),f^X(y))\mid (z,x,y) \in S^X\}.]$\pars

For the infinite hyperfinite model, $M$ would need to be, at least, denumerable. This could be accomplished using systems that vary due to parameter changes. The set $\hyper {M} - M$ 
forms a higher-language and members can have significant partial meanings (Herrmann, 1993, p. 101-102). Theoretically, for any injection $f^H\colon [0,\lambda +1] \to \hyper M,$ the structure $f^H[{\cal D}_\lambda^H]$ models $TAS_1$ and $TAS_2$.\pars

Most certainly these models do not constitute all of specific structures that model each scheme. For example, Gagnon's model for each $TAS_i,\ i = 1,2,3$ is distinct from the ones presented here. A particular dialectical argument may require entirely different structure relations than any presented in this paper. \parm

\centerline{\bf References}\parm 
Gagnon, L.S. 1980. {\it Three theories of dialectic}, Notre Dame Journal of 
Formal Logic, XXI (2):316--318.\pars
Herrmann, R. A. 2006a. {\it General logic-systems and finite consequence operators}, Logica Universalis, 1:201-208. http://arxiv.org/abs/math/0512559\pars
Herrmann, R. A. 2006b. The GGU-model and generation of developmental paradigms, http//:arxiv.org/abs/math/0605120 \pars
Herrmann, R. A. 2001. {\it Hyperfinite and standard unifications for physical theories}, Internat. J. Math. and Math. Sci., 28(2):93-102. http://arxiv.org/abs/physics/0105012 \pars
Herrmann, R. A. 1993. The Theory of Ultralogics.\hfil\break http://arxiv.org/abs/math/9903081 http://arxiv.org/abs/math/9903082\par 
Herrmann, R. A. 1991. Nonstandard Analysis Applied to Advanced Undergraduate Mathematics - Infinitesimal Modeling. http://arxiv.org/abs/math/03122432\pars

\bigskip

\noindent Robert A. Herrmann, USNA, 43017 Daventry Sq, South Riding, VA, 20152-2050, USA\pars
\noindent {\it E-mail address:} herrmann@serve.com \end
\end